\documentclass[11pt]{article}
\usepackage{amsfonts}
\usepackage{mathrsfs}
\usepackage[dvips]{graphics}
\usepackage[dvips]{color}
\usepackage{amsthm}
\usepackage[]{amsmath}
\usepackage{CJK}
\usepackage{indentfirst}
\newtheorem{proposition}{Proposition}[section]

\newtheorem{theorem}[proposition]{Theorem}
\newtheorem*{tA}{Theorem A}
\newtheorem*{tB}{Theorem B}
\newtheorem*{tC}{Theorem C}
\newtheorem*{tD}{Theorem D}
\newtheorem*{aB}{Algorithm B}
\newtheorem*{aD}{Algorithm D}

\newtheorem{property}[proposition]{Property}

\begin{document}
\begin{CJK*}{GBK}{song}
\CJKindent

\centerline{\textbf{\LARGE{The numbers of distinct and repeated}}}

\vspace{0.2cm}

\centerline{\textbf{\LARGE{squares and cubes in the Tribonacci sequence}}}

\vspace{0.2cm}

\centerline{Huang Yuke\footnote[1]{School of Mathematics and Systems Science, Beihang University (BUAA), Beijing, 100191, P. R. China. E-mail address: huangyuke@buaa.edu.cn,~hyg03ster@163.com (Corresponding author).}
~~Wen Zhiying\footnote[2]{Department of Mathematical Sciences, Tsinghua University, Beijing, 100084, P. R. China. E-mail address: wenzy@tsinghua.edu.cn.}}

\vspace{1cm}

\noindent\textbf{Abstract:} The Tribonacci sequence $\mathbb{T}$ is the fixed point of the substitution $\sigma(a,b,c)=(ab,ac,a)$.
The main result is twofold: (1)
we give the explicit expressions of the numbers of distinct squares and cubes in $\mathbb{T}[1,n]$ (the prefix of $\mathbb{T}$ of length $n$); (2) we give algorithms for counting the number of repeated squares and cubes in $\mathbb{T}[1,n]$ for all $n$; then get explicit expressions for some special $n$ such as $n=t_m$ (the Tribonacci number).

\vspace{0.2cm}


\vspace{0.2cm}

\noindent\textbf{2010 MR Subject Classification}~~~~11B85; 68Q45

\section{Introduction}

Let $\mathcal{A}=\{a,b,c\}$ be a three-letter alphabet.
The Tribonacci sequence $\mathbb{T}$ is the fixed point of the substitution $\sigma(a,b,c)=(ab,ac,a)$.
As a natural generalization of the Fibonacci sequence, the Tribonacci sequence has
been studied extensively by many authors, see \cite{HW2015-2,MS2014,TW2007}.

Let $\omega$ be a factor of $\mathbb{T}$, denoted by $\omega\prec\mathbb{T}$.
Let $\omega_p$ be the $p$-th occurrence of $\omega$.
If the factor $\omega$ and integer $p$ such that $\omega_p\omega_{p+1}\prec\mathbb{T}$,
we call $\omega_p\omega_{p+1}$ a square of $\mathbb{T}$.
We call $\omega\omega$ a square of $\mathbb{T}$ if there exist $p$ such that $\omega_p\omega_{p+1}\prec\mathbb{T}$.
Similarly, we call $\omega_p\omega_{p+1}\omega_{p+2}$ (resp. $\omega\omega\omega$) a cube of $\mathbb{T}$.
As we know, $\mathbb{T}$ contains no fourth powers. The properties of squares and cubes are objects of a great interest in many aspects of mathematics and computer science etc.

We denote by $|\omega|$ the length of $\omega$, and by $|\omega|_\alpha$ the number of letter $\alpha$ in $\omega$, $\alpha\in\mathcal{A}$.
Denote $P(\omega,p)$ the position of the last letter of $\omega_p$.
Let $\tau=x_1\cdots x_n$. % be a finite word (or $\tau=x_1x_2\cdots$ be a sequence).
For $i\leq j\leq n$, define $\tau[i,j]=x_ix_{i+1}\cdots x_{j-1}x_j$.
By convention, denote $\tau[i]=\tau[i,i]=x_i$, $\tau[i,i-1]=\varepsilon$ (empty word).
Denote $T_m=\sigma^m(a)$ for $m\geq0$, $T_{-2}=\varepsilon$, $T_{-1}=c$,
then $T_0=a$, $T_1=ab$ and $T_{m}=T_{m-1}T_{m-2}T_{m-3}$ for $m\geq2$.
Denote $t_m=|T_m|$ for $m\geq-2$, called the $m$-th Tribonacci number.
%Then $t_{-2}=0$, $t_{-1}=t_{0}=1$, $t_m=t_{m-1}+t_{m-2}+t_{m-3}$ for $m\geq1$.
Denote by $\delta_m$ the last letter of $T_m$, then $\delta_m=a$ (resp. $b$, $c$) for $m\equiv0$, (resp. $1$, $2$) mod 3 and $m\geq-1$.

Let $\mathbb{T}[1,n]$ be the prefix of $\mathbb{T}$ of length $n$. In this paper, we consider the four functions below:

$A(n)=\sharp\{\omega:\omega\omega\prec\mathbb{T}[1,n]\}$, the number of distinct squares in $\mathbb{T}[1,n]$;

$B(n)=\sharp\{(\omega,p):\omega_p\omega_{p+1}\prec\mathbb{T}[1,n]\}$, the number of repeated squares in $\mathbb{T}[1,n]$;

$C(n)=\sharp\{\omega:\omega\omega\omega\prec\mathbb{T}[1,n]\}$, the number of distinct cubes in $\mathbb{T}[1,n]$;

$D(n)=\sharp\{(\omega,p):\omega_p\omega_{p+1}\omega_{p+2}\prec\mathbb{T}[1,n]\}$, the number of repeated cubes in $\mathbb{T}[1,n]$.

In 2006, A.Glen\cite{G2006} gave expressions of $A(t_m)$.
In 2014, H.Mousavi and J.Shallit\cite{MS2014} gave expressions of $B(t_m)$ and $D(t_m)$, which they proved by mechanical way.
All of these results above only consider the squares or cubes in the prefixes of some special lengths: the Tribonacci numbers.
In this paper, for all $n,m\geq1$, we give: (1)
explicit expressions of $A(n)$, $C(n)$, $B(t_m)$ and $D(t_m)$; (2) fast algorithms for counting $B(n)$ and $D(n)$.
We also considered this problem in the Fibonacci sequence, see \cite{HW2016-2}.
%\section{Main Results}

We denote the gap between $\omega_p$ and $\omega_{p+1}$ by $G_p(\omega)$.
The sequence $\{G_p(\omega)\}_{p\ge 1}$ is called the gap sequence of factor $\omega$.
Taking ``square" for example, the \textbf{main difficulty} is twofold:
(1) The positions of all squares are not easy to be determined.
We overcome this difficulty by using the ``gap sequence" property of $\mathbb{T}$, which we introduced and studied in \cite{HW2015-2}.
(2) By the gap sequence property of $\mathbb{T}$, we can find out all distinct squares in $\mathbb{T}[1,n]$. We can also count the number of occurrences of each square. So the summation of these numbers are the numbers of repeated squares in $\mathbb{T}[1,n]$. But this method is complicated.
We overcome this difficulty by studying the relations among positions of each $\omega_p$, and establishing a recursive structure, called square trees.
First we list the main results as below.

\begin{tA}[The numbers of distinct squares, $A(n)$]\label{A}\
$A(n)=0$ for $n\leq7$; $A(n)=1$ for $n=8,9$; $A(n)=2$ for $10\leq n\leq13$.
For $n\geq14$, let $m$ such that $\alpha_m\leq n<\alpha_{m+1}$, then $m\geq4$,
\begin{equation*}
A(n)=
\begin{cases}
n-\frac{1}{2}(t_{m}+t_{m-3}+m+3),&\alpha_m\leq n<\beta_{m};\\
\frac{1}{2}(t_{m-1}+t_{m-2}+4t_{m-3}-m-5),&\beta_{m}\leq n<\gamma_m;\\
n-\frac{1}{2}(t_{m-1}+3t_{m-2}+m+3),&\gamma_m\leq n<\theta_m;\\
\frac{1}{2}(2t_{m-1}+t_{m-2}+3t_{m-3}-m-6),&\theta_m\leq n<\alpha_{m+1}.
\end{cases}
\end{equation*}
Here $\alpha_m=2t_{m-1}$, $\beta_m=t_{m}+2t_{m-3}-1$, $\gamma_m=2t_{m}-t_{m-1}$ and $\theta_m=\frac{3t_{m}+t_{m-2}-3}{2}$.
\end{tA}

\noindent\emph{Example.} Consider $n=65\in[\alpha_6,\cdots,\alpha_{7}-1]=[48,\cdots,87]$, $m=6$.
Moreover $\gamma_6=2t_{6}-t_{5}=64\leq n<\theta_6=\frac{3t_{6}+t_{4}-3}{2}=71$,
$A(65)=65-\frac{1}{2}(t_{5}+3t_{4}+6+3)=29$.
In fact, the positions of the last letters of the 29 squares are
$\{8,10,14,15,16,19,20,\underbrace{26,\cdots,31}_6,\underbrace{35,\cdots,38}_4,
\underbrace{48,\cdots,57}_{10},64,65\}.$
Here the number under ``$\stackrel{\underbrace{}}{}$" means the number of elements.

\begin{tB}[The numbers of repeated squares, $B(t_m)$]\label{B} For $m\geq3$, $$B(t_{m})=\tfrac{m}{22}(9t_m-t_{m-1}-5t_{m-2})
+\tfrac{1}{44}(-81t_m+26t_{m-1}+13t_{m-2})+m+\tfrac{1}{4}.$$
\end{tB}

\noindent\emph{Example.} $B(t_{5})=B(24)=\tfrac{5}{22}(9t_5-t_{4}-5t_{3})
+\tfrac{1}{44}(-81t_5+26t_{4}+13t_{3})+5+\tfrac{1}{4}=9$.

\begin{tC}[The numbers of distinct cubes, $C(n)$]\label{C}\
$C(n)=0$ for $n\leq57$.
For $n\geq58$, let $m$ such that $t_{m-1}+2t_{m-4}\leq n<t_{m}+2t_{m-3}$, then $m\geq7$,
\begin{equation*}
C(n)=
\begin{cases}
n-\frac{1}{2}(4t_{m-1}-t_{m-2}-3t_{m-3}+m-6),&n\leq\frac{3t_{m-1}-t_{m-3}-3}{2};\\
\frac{1}{2}(t_{m-5}+t_{m-6}-m+3),&otherwise.
\end{cases}
\end{equation*}
\end{tC}

\noindent\emph{Example.} Consider $n=365\in[t_{9}+2t_{6},\cdots,t_{10}+2t_{7}]=[362,\cdots,666]$, $m=10$.
Moreover $n\leq\tfrac{3t_{9}-t_{7}-3}{2}=369$,
$C(365)=365-\frac{1}{2}(4t_{9}-t_{8}-3t_{7}+10-6)=11$.
In fact, the positions of the last letters of the 11 cubes are
$\{58,107,108,197,198,199,200,362,363,364,365\}.$

\begin{tD}[The numbers of repeated cubes, $D(t_m)$]\label{D} For $m\geq3$,
$$\begin{array}{rl}
D(t_m)=&\tfrac{m}{22}(-6t_m+8t_{m-1}+7t_{m-2})
+\tfrac{1}{44}(-23t_m+34t_{m-1}-5t_{m-2})+\tfrac{m}{6}\\
&~~~~-\frac{1}{4}[m\equiv0(\mathrm{mod}~3)]+\frac{1}{12}[m\equiv1(\mathrm{mod}~3)]
+\frac{5}{12}[m\equiv2(\mathrm{mod}~3)].
\end{array}$$
Here $[P]$ is Iverson notation, and equals 1 if P holds and 0 otherwise.
\end{tD}

\noindent\emph{Example.} $D(t_8)=D(149)=\tfrac{8}{22}(-6t_8+8t_{7}+7t_{6})
+\tfrac{1}{44}(-23t_8+34t_{7}-5t_{6})+\tfrac{8}{6}+\frac{5}{12}=4$.

\vspace{0.2cm}

H.Mousavi and J.Shallit gave \textbf{Theorem B.} and \textbf{Theorem D.} in \cite{MS2014}.
We prove them as corollaries of  \textbf{Algorithm B.} and \textbf{Algorithm D.}, which
are fast algorithms for counting $B(n)$ and $D(n)$ for all $n$, see Section 4 and 6, respectively.

\section{Preliminaries of the Tribonacci sequence}

We define the kernel numbers that $k_{0}=0$, $k_{1}=k_{2}=1$, $k_m=k_{m-1}+k_{m-2}+k_{m-3}-1$ for $m\geq3$.
Then $k_m=k_{m-3}+t_{m-4}=\frac{t_{m-3}+t_{m-5}+1}{2}$ for $m\geq3$.
The kernel word with order $m$ is defined as
$K_1=a$, $K_2=b$, $K_3=c$, $K_m=\delta_{m-1}T_{m-3}[1,k_m-1]$ for $m\geq4$.
Let $Ker(\omega)$ be the maximal kernel word occurring in factor $\omega$, then by Theorem 4.3 in \cite{HW2015-2}, $Ker(\omega)$ occurs in $\omega$ only once.
Moreover,

\begin{property}[Theorem 4.11 in \cite{HW2015-2}]\label{wpt}
$Ker(\omega_p)=Ker(\omega)_p$ for all $\omega\in\mathbb{T}$ and $p\geq1$.
\end{property}

\begin{property}[Theorem 3.3 in \cite{HW2015-2}]\label{Gt}
The gap sequence $\{G_p(\omega)\}_{p\geq1}$ is the Tribonacci sequence over the alphabet $\{G_1(\omega)$,$G_2(\omega)$,$G_4(\omega)\}$ for all $\omega\in\mathbb{T}$.
\end{property}

\begin{property}[Property 6.1 in \cite{HW2016-3}]\label{P} For $m,p\geq1$,

$P(K_m,p)=pt_{m-1}+|\mathbb{T}[1,p-1]|_a(t_{m-2}+t_{m-3})+|\mathbb{T}[1,p-1]|_bt_{m-2}+k_m-1$.
\end{property}

\noindent\emph{Example.}
$P(K_m,1)=k_{m+3}-1=\frac{t_{m}+t_{m-2}-1}{2}$ for $m\geq1$.
$P(a,p)=p+|\mathbb{T}[1,p-1]|_a+|\mathbb{T}[1,p-1]|_b$,
$P(b,p)=2p+2|\mathbb{T}[1,p-1]|_a+|\mathbb{T}[1,p-1]|_b$,
$P(c,p)=4p+3|\mathbb{T}[1,p-1]|_a+2|\mathbb{T}[1,p-1]|_b$ for $p\geq1$.

\begin{property}[Lemma 6.4. in \cite{HW2016-3}]\label{L} For $p\geq1$,
(1) $|\mathbb{T}[1,P(\alpha,p)]|_\alpha=p$ for $\alpha\in\mathcal{A}$;

(2) $|\mathbb{T}[1,P(b,p)]|_a=P(a,p)$, $|\mathbb{T}[1,P(a,p)|_b=|\mathbb{T}[1,p-1]|_a$;

(3) $|\mathbb{T}[1,P(c,p)]|_a=P(b,p)$, $|\mathbb{T}[1,P(c,p)]|_b=P(a,p)$.
\end{property}

\section{The number of distinct squares, $A(n)$}

By Lemma 4.7, Definition 4.12 and Corollary 4.13 in \cite{HW2015-2}, any factor $\omega$ with kernel $K_m$ can be expressed uniquely as
$\omega=T_{m-1}[i,t_{m-1}-1]K_mT_{m}[k_m,k_m+j-1],$
where $1\leq i\leq t_{m-1}$ and $0\leq j\leq t_{m-1}-1$.

By Theorem 3.3, Corollary 3.12 and Proposition 6.7(1) in \cite{HW2015-2}, $\omega_p\omega_{p+1}\prec\mathbb{T}$ has three cases.

\vspace{0.2cm}

\textbf{Case 1.} $G_p(K_m)=G_1(K_m)$. Here $|G_1(K_m)|=t_{m}-k_{m}$.

Since $|G_p(K_m)|=|T_{m-1}[i,t_{m-1}-1]|+|T_{m}[k_m,k_m+j-1]|=t_{m-1}-i+j$, $j=t_{m-2}+t_{m-3}-k_{m}+i$. So $0\leq j\leq t_{m-1}-1$ gives a range of $i$. Comparing this range with $1\leq i\leq t_{m-1}$, we have $1\leq i\leq k_{m+1}-1$ and $m\geq3$.
Furthermore,
\begin{equation*}
\begin{split}
\omega=&T_{m-1}[i,t_{m-1}]T_{m}[1,t_{m-2}+t_{m-3}+i-1]
=T_{m-1}[i,t_{m-1}]T_{m-2}T_{m-3}T_{m-2}[1,i-1];\\
\omega\omega
=&T_{m}[i,t_{m}-1]\underline{\delta_{m}T_{m-1}[1,k_{m+1}-1]}T_{m+1}[k_{m+1},t_{m}+i-1].
\end{split}
\end{equation*}
Thus $K_{m+1}=\delta_{m}T_{m-1}[1,k_{m+1}-1]\prec\omega\omega$.
Similarly, $K_{m+2},K_{m+3},K_{m+4}\!\not\prec\omega\omega$ and $|\omega\omega|<K_{m+5}$, so $K_{m+1}$ is the largest kernel word in $\omega\omega$, i.e. $Ker(\omega\omega)=K_{m+1}$ for $m\geq3$.
Moreover, since $|\omega|=|G_p(K_m)|+k_m$, $|\omega|=t_{m}$.

\vspace{0.2cm}

By analogous arguments, we have

\textbf{Case 2.} $G_p(K_m)=G_2(K_m)$. Here $|G_2(K_m)|=t_{m-2}+t_{m-1}-k_{m}$.
$$\omega\omega
=T_{m}[i,t_{m-1}+t_{m-2}-1]\underline{K_{m+2}}T_{m+1}[k_{m+2},t_{m-1}+t_{m-2}+i-1],$$
where $1\leq i\leq k_{m+2}-1$, $m\geq2$, $Ker(\omega\omega)=K_{m+2}$ and $|\omega|=t_{m-2}+t_{m-1}$.

\textbf{Case 3.} $G_p(K_m)=G_4(K_m)$. Here $|G_4(K_m)|=t_{m-1}-k_{m}$.
$$\omega\omega
=T_{m-1}[i,t_{m-1}-1]\underline{K_{m+3}}T_{m+1}[k_{m+3},t_{m-1}+i-1],$$
where $k_{m}\leq i\leq t_{m-1}$, $m\geq1$, $Ker(\omega\omega)=K_{m+3}$ and $|\omega|=t_{m-1}$.

\vspace{0.2cm}

\noindent\emph{Remark.}
By the three cases of squares, we have: (1) all squares in $\mathbb{T}$ are of length $2t_m$ or $2t_m+2t_{m-1}$ for some $m\geq0$; (2) for all $m\geq0$, there exists a square of length $2t_m$ and $2t_m+2t_{m-1}$ in $\mathbb{T}$. These are known results of H.Mousavi and J.Shallit, see Theorem 5 in \cite{MS2014}.

\vspace{0.2cm}

We define three sets for $m\geq4$,
\begin{equation*}
\begin{cases}
\langle1,K_m,p\rangle&=
\{P(\omega\omega,p):Ker(\omega\omega)=K_m,|\omega|=t_{m-1},\omega\omega\prec\mathbb{T}\}\\
&=\{P(K_m,p)+\frac{-t_{m}+4t_{m-1}-t_{m-2}+1}{2},\cdots,P(K_m,p)+t_{m-1}-1\};\\
\langle2,K_m,p\rangle&=
\{P(\omega\omega,p):Ker(\omega\omega)=K_m,|\omega|=t_{m-4}+t_{m-3},\omega\omega\prec\mathbb{T}\}\\
&=\{P(K_m,p)+\frac{-t_{m}+4t_{m-1}-3t_{m-2}+1}{2},\cdots,
P(K_m,p)+t_{m-1}-t_{m-2}-1\};\\
\langle3,K_m,p\rangle&=
\{P(\omega\omega,p):Ker(\omega\omega)=K_m,|\omega|=t_{m-4},\omega\omega\prec\mathbb{T}\}\\
&=\{P(K_m,p),\cdots,P(K_m,p)+\frac{-5t_{m}+10t_{m-1}-t_{m-2}-1}{2}\}.
\end{cases}
\end{equation*}
Obviously, these sets correspond the positions $P(\omega\omega,p)$ for the three cases of squares, respectively. Each set contains some consecutive integers.
Moreover $\sharp\langle 1,K_m,p\rangle=\sharp\langle 2,K_m,p\rangle=\frac{t_{m}-2t_{m-1}+t_{m-2}-1}{2}$,
$\sharp\langle 3,K_m,p\rangle=\frac{-5t_{m}+10t_{m-1}-t_{m-2}+1}{2}$.

Since $P(K_m,1)=\frac{t_{m}+t_{m-2}-1}{2}$, we have
$\langle1,K_m,1\rangle=\{2t_{m-1},\cdots,\frac{t_{m}+2t_{m-1}+t_{m-2}-3}{2}\}$,
$\langle2,K_m,1\rangle=\{2t_{m-1}-t_{m-2},\cdots,\frac{t_{m}+2t_{m-1}-t_{m-2}-3}{2}\}$ and
$\langle3,K_m,1\rangle=\{\frac{t_{m}+t_{m-2}-1}{2},\cdots,-2t_{m}+5t_{m-1}-1\}$.
Obviously, $\langle j,K_m,1\rangle$ are pairwise disjoint.
Moreover, since $\max\langle1,K_m,1\rangle+1=\min\langle3,K_{m+1},1\rangle$, sets $\langle1,K_m,1\rangle$ and $\langle3,K_{m+1},1\rangle$ are consecutive.
We have
$\langle1,K_m,1\rangle\cup\langle3,K_{m+1},1\rangle=\{2t_{m-1},\cdots,t_{m}+2t_{m-3}-1\}$.

Therefore we get a chain
$$\langle3,K_4,1\rangle,\langle2,K_4,1\rangle,\langle1,K_4,1\rangle\cup\langle3,K_5,1\rangle,\cdots,
\langle2,K_m,1\rangle,\langle1,K_m,1\rangle\cup\langle3,K_{m+1},1\rangle,\cdots$$

Denote $a(n)=\sharp\{\omega:\omega\omega\triangleright\mathbb{T}[1,n],
\omega\omega\not\!\prec\mathbb{T}[1,n-1]\}$, then $A(n)=\sum_{i=1}^n a(i)$.

\begin{property}[]\label{a}~For $n<14$, $a(n)=1$ if and only if $n\in\{8,10\}$;
for $n\geq14$, let $m$ such that $2t_{m-1}\leq n<2t_{m}$, then $m\geq4$ and $a(n)=1$ if and only if
$$n\in\{2t_{m-1},\cdots,t_{m}+2t_{m-3}-1\}
\cup\{2t_{m}-t_{m-1},\cdots,\tfrac{3t_{m}+t_{m-2}-3}{2}\}.$$
\end{property}

For $m\geq3$, we denote
$\Delta_m=\sum_{i=4}^m\sharp\langle j,K_i,1\rangle$ for $j=1,2$, $\Theta_m=\sum_{i=4}^m\sharp\langle 3,K_i,1\rangle$.
So $\Delta_3=\Theta_3=0$.
And for $m\geq4$, since $\sum_{i=0}^{m}t_i=\frac{t_m+t_{m+2}-3}{2}$ and $\sum_{i=1}^{m}k_i=\frac{t_{m-2}+t_{m-3}+m}{2}$,
we have $\Delta_{m}=\frac{t_{m-2}+t_{m-3}-m}{2}$ and $\Theta_m=\frac{t_{m-2}-t_{m-3}+2t_{m-4}+m-6}{2}$.
For $m\geq4$, we denote
\begin{equation*}
\begin{cases}
\alpha_m=\min\langle1,K_m,1\rangle=2t_{m-1},~\beta_m=\max\langle3,K_{m+1},1\rangle=t_{m}+2t_{m-3}-1,\\
\gamma_m=\min\langle2,K_{m+1},1\rangle=2t_{m}-t_{m-1},
\theta_m=\max\langle2,K_{m+1},1\rangle=\frac{3t_{m}+t_{m-2}-3}{2}.
\end{cases}
\end{equation*}

By Property \ref{a} and the definition of $\Delta_m$, $\Theta_m$,
\begin{equation*}
\begin{cases}
A(\alpha_m)=\Delta_{m-1}+\Delta_{m}+\Theta_{m}+1
=\frac{2t_{m-2}+t_{m-3}+3t_{m-4}-m-3}{2},\\
A(\beta_m)=\Delta_{m}+\Delta_{m}+\Theta_{m+1}
=\frac{t_{m-1}+t_{m-2}+4t_{m-3}-m-5}{2},\\
A(\gamma_m)=\Delta_{m}+\Delta_{m}+\Theta_{m+1}+1
=A(\beta_m)+1,\\
A(\theta_m)=\Delta_{m}+\Delta_{m+1}+\Theta_{m+1}
=A(\alpha_{m+1})-1.
\end{cases}
\end{equation*}
When $\alpha_m\leq n<\beta_m$, $A(n)=A(\alpha_m)+n-\alpha_m$;
when $\beta_m\leq n<\gamma_m$, $A(n)=A(\beta_m)$;
when $\gamma_m\leq n<\theta_m$, $A(n)=A(\gamma_m)+n-\gamma_m$;
when $\theta_m\leq n<\alpha_{m+1}$, $A(n)=A(\theta_m)$. So
we get \textbf{Theorem A}.

\vspace{0.2cm}

For $m\geq5$, since $\theta_{m-1}=\frac{3t_{m-1}+t_{m-3}-1}{2}\leq t_m<\alpha_{m}$, $A(t_m)=A(\theta_{m-1})$ by Theorem A.
It is easy to check the expression holds also for $m=3,4$. Thus for $m=0,1,2$, $A(t_m)=0$, and

\begin{theorem}[]\
$A(t_m)=\frac{1}{2}(2t_{m-2}+t_{m-3}+3t_{m-4}-m-5)$ for $m\geq3$.
\end{theorem}

\noindent\emph{Remark.} A.Glen gave $A(t_m)$ in Theorem 6.30 in \cite{G2006}, that
$A(t_m)=\sum_{i=0}^{m-2}(d_i+1)+d_{m-4}+d_{m-5}+1$ for $m\geq3$,
where $d_{-2}=d_{-1}=-1$, $d_0=0$ and $d_m=\frac{t_{m+1}+t_{m-1}-3}{2}$ for $m\geq1$.
Since $\sum_{i=0}^{m}t_i=\frac{t_m+t_{m+2}-3}{2}$, the two expressions of $A(t_m)$ are same.

\section{The number of repeated squares, $B(n)$}

For $m\geq4$ and $p\geq1$, we consider the vectors
\begin{equation*}
\begin{cases}
\Gamma_{1,m,p}=[P(K_m,p)+t_{m-1}-t_{m-2},\cdots,P(K_m,p)+t_{m-1}-1];\\
\Gamma_{2,m,p}=[P(K_m,p)-t_{m}+2t_{m-1},\cdots,P(K_m,p)+t_{m-1}-t_{m-2}-1];\\
\Gamma_{3,m,p}=[P(K_m,p),\cdots,P(K_m,p)-t_{m}+2t_{m-1}-1].
\end{cases}
\end{equation*}
Obviously, $\max\Gamma_{2,m,p}+1=\min\Gamma_{1,m,p}$ and
$\max\Gamma_{3,m,p}+1=\min\Gamma_{2,m,p}$.
Moreover $|\Gamma_{1,m,p}|=t_{m-2}$, $|\Gamma_{2,m,p}|=t_{m-3}$ and $|\Gamma_{3,m,p}|=t_{m-4}$.
Using Property \ref{L},
comparing minimal and maximal elements in these sets below, we have
\begin{equation*}
\begin{cases}
\Gamma_{1,m,p}=[\Gamma_{3,m-1,P(a,p)+1},\Gamma_{2,m-1,P(a,p)+1},\Gamma_{1,m-1,P(a,p)+1}],
&m\geq5;\\
\Gamma_{2,m,p}=[\Gamma_{3,m-2,P(b,p)+1},\Gamma_{2,m-2,P(b,p)+1},\Gamma_{1,m-2,P(b,p)+1}],
&m\geq6;\\
\Gamma_{3,m,p}=[\Gamma_{3,m-3,P(c,p)+1},\Gamma_{2,m-3,P(c,p)+1},\Gamma_{1,m-3,P(c,p)+1}],
&m\geq7.
\end{cases}
\end{equation*}
Thus we establish recursive relations for any $\Gamma_{1,m,p}$ ($m\geq5$),
$\Gamma_{2,m,p}$ ($m\geq6$) and $\Gamma_{3,m,p}$ ($m\geq7$).

It is easy to check that $\langle 1,K_m,p\rangle$ (resp. $\langle 2,K_m,p\rangle$, $\langle 3,K_m,p\rangle$) contains the several maximal (resp. maximal, minimal) elements of $\Gamma_{1,m,p}$ (resp. $\Gamma_{2,m,p}$, $\Gamma_{3,m,p}$).
Thus we get recursive relations of the positions of repeated squares in $\mathbb{T}$, called square trees.
Here we denote $P(\alpha,p)+1$ by $\hat{\alpha}$ for short, $\alpha\in\{a,b,c\}$.
For $i\in\{1,2,3\}$, $\pi_i$ is a substitution over $\{\langle j,K_m,p\rangle:m\geq4+j,p\geq1\}$ that
\begin{equation*}
\begin{cases}
\pi_1\langle1,K_m,p\rangle=\langle3,K_{m-1},\hat{a}\rangle\cup\langle2,K_{m-1},\hat{a}\rangle
\cup\langle1,K_{m-1},\hat{a}\rangle,&m\geq5;\\
\pi_2\langle2,K_m,p\rangle=\langle3,K_{m-2},\hat{b}\rangle\cup\langle2,K_{m-2},\hat{b}\rangle
\cup\langle1,K_{m-2},\hat{b}\rangle,&m\geq6;\\
\pi_3\langle3,K_m,p\rangle=\langle3,K_{m-3},\hat{c}\rangle\cup\langle2,K_{m-3},\hat{c}\rangle
\cup\langle1,K_{m-3},\hat{c}\rangle,&m\geq7.
\end{cases}
\end{equation*}

\begin{figure}
\footnotesize
\setlength{\unitlength}{1mm}
\begin{center}
\begin{picture}(140,43)
\put(6,1){51}
\put(6,4){50}
\put(6,7){49}
\put(6,10){48}
\put(1,13){$\langle1,K_6,1\rangle$}
\put(1,0){\line(0,1){16}}
\put(1,0){\line(1,0){13}}
\put(14,0){\line(0,1){16}}
\put(1,16){\line(1,0){13}}
\put(26,1){51}
\put(26,4){50}
\put(21,7){$\langle1,K_5,2\rangle$}
\put(21,0){\line(0,1){10}}
\put(21,0){\line(1,0){13}}
\put(34,0){\line(0,1){10}}
\put(21,10){\line(1,0){13}}
\put(46,22){44}
\put(46,25){43}
\put(41,28){$\langle2,K_5,2\rangle$}
\put(41,21){\line(0,1){10}}
\put(41,21){\line(1,0){13}}
\put(54,21){\line(0,1){10}}
\put(41,31){\line(1,0){13}}
\put(66,34){40}
\put(66,37){39}
\put(61,40){$\langle3,K_5,2\rangle$}
\put(61,33){\line(0,1){10}}
\put(61,33){\line(1,0){13}}
\put(74,33){\line(0,1){10}}
\put(61,43){\line(1,0){13}}
\put(86,1){51}
\put(81,4){$\langle1,K_4,4\rangle$}
\put(81,0){\line(0,1){7}}
\put(81,0){\line(1,0){13}}
\put(94,0){\line(0,1){7}}
\put(81,7){\line(1,0){13}}
\put(106,13){47}
\put(101,16){$\langle2,K_4,4\rangle$}
\put(101,12){\line(0,1){7}}
\put(101,12){\line(1,0){13}}
\put(114,12){\line(0,1){7}}
\put(101,19){\line(1,0){13}}
\put(126,19){45}
\put(121,22){$\langle3,K_4,4\rangle$}
\put(121,18){\line(0,1){7}}
\put(121,18){\line(1,0){13}}
\put(134,18){\line(0,1){7}}
\put(121,25){\line(1,0){13}}
\put(15,10){\vector(1,-1){5}}
\put(15,10){\vector(2,1){25}}
\put(15,10){\line(1,1){25}}
\put(40,35){\vector(1,0){20}}
\put(16,15){$\pi_1$}
\put(38,7){$\pi_1$}
\put(35,4){\vector(1,0){45}}
\put(35,4){\line(5,1){50}}
\put(85,14){\vector(1,0){15}}
\put(35,4){\line(3,1){54}}
\put(89,22){\vector(1,0){31}}
\put(120,0){\footnotesize{\textbf{(a)}}}
\end{picture}
\begin{picture}(140,47)
\put(6,1){71}
\put(6,4){70}
\put(6,7){69}
\put(6,10){68}
\put(6,13){67}
\put(6,16){66}
\put(6,19){65}
\put(6,22){64}
\put(1,25){$\langle2,K_7,1\rangle$}
\put(1,0){\line(0,1){28}}
\put(1,0){\line(1,0){13}}
\put(14,0){\line(0,1){28}}
\put(1,28){\line(1,0){13}}
\put(26,1){71}
\put(26,4){70}
\put(21,7){$\langle1,K_5,3\rangle$}
\put(21,0){\line(0,1){10}}
\put(21,0){\line(1,0){13}}
\put(34,0){\line(0,1){10}}
\put(21,10){\line(1,0){13}}
\put(46,22){64}
\put(46,25){63}
\put(41,28){$\langle2,K_5,3\rangle$}
\put(41,21){\line(0,1){10}}
\put(41,21){\line(1,0){13}}
\put(54,21){\line(0,1){10}}
\put(41,31){\line(1,0){13}}
\put(66,34){60}
\put(66,37){59}
\put(61,40){$\langle3,K_5,3\rangle$}
\put(61,33){\line(0,1){10}}
\put(61,33){\line(1,0){13}}
\put(74,33){\line(0,1){10}}
\put(61,43){\line(1,0){13}}
\put(86,1){71}
\put(81,4){$\langle1,K_4,6\rangle$}
\put(81,0){\line(0,1){7}}
\put(81,0){\line(1,0){13}}
\put(94,0){\line(0,1){7}}
\put(81,7){\line(1,0){13}}
\put(106,13){67}
\put(101,16){$\langle2,K_4,6\rangle$}
\put(101,12){\line(0,1){7}}
\put(101,12){\line(1,0){13}}
\put(114,12){\line(0,1){7}}
\put(101,19){\line(1,0){13}}
\put(126,19){65}
\put(121,22){$\langle3,K_4,6\rangle$}
\put(121,18){\line(0,1){7}}
\put(121,18){\line(1,0){13}}
\put(134,18){\line(0,1){7}}
\put(121,25){\line(1,0){13}}
\put(15,10){\vector(1,-1){5}}
\put(15,10){\vector(2,1){25}}
\put(15,10){\line(1,1){25}}
\put(40,35){\vector(1,0){20}}
\put(16,15){$\pi_2$}
\put(38,7){$\pi_1$}
\put(35,4){\vector(1,0){45}}
\put(35,4){\line(5,1){50}}
\put(85,14){\vector(1,0){15}}
\put(35,4){\line(3,1){54}}
\put(89,22){\vector(1,0){31}}
\put(120,0){\footnotesize{\textbf{(b)}}}
\end{picture}
\begin{picture}(140,47)
\put(5,7){106}
\put(5,10){105}
\put(5,13){104}
\put(5,16){103}
\put(5,19){102}
\put(5,22){101}
\put(5,25){100}
\put(6,28){99}
\put(6,31){98}
\put(6,34){97}
\put(6,37){96}
\put(1.5,40){$\langle3,K_8,1\rangle$}
\put(1,6){\line(0,1){37}}
\put(1,6){\line(1,0){13}}
\put(14,6){\line(0,1){37}}
\put(1,43){\line(1,0){13}}
\put(25,1){108}
\put(25,4){107}
\put(21.5,7){$\langle1,K_5,5\rangle$}
\put(21,0){\line(0,1){10}}
\put(21,0){\line(1,0){13}}
\put(34,0){\line(0,1){10}}
\put(21,10){\line(1,0){13}}
\put(45,22){101}
\put(45,25){100}
\put(41.5,28){$\langle2,K_5,5\rangle$}
\put(41,21){\line(0,1){10}}
\put(41,21){\line(1,0){13}}
\put(54,21){\line(0,1){10}}
\put(41,31){\line(1,0){13}}
\put(66,34){97}
\put(66,37){96}
\put(61.5,40){$\langle3,K_5,5\rangle$}
\put(61,33){\line(0,1){10}}
\put(61,33){\line(1,0){13}}
\put(74,33){\line(0,1){10}}
\put(61,43){\line(1,0){13}}
\put(85,1){108}
\put(81.5,4){$\langle1,K_4,9\rangle$}
\put(81,0){\line(0,1){7}}
\put(81,0){\line(1,0){13}}
\put(94,0){\line(0,1){7}}
\put(81,7){\line(1,0){13}}
\put(105,13){104}
\put(101.5,16){$\langle2,K_4,9\rangle$}
\put(101,12){\line(0,1){7}}
\put(101,12){\line(1,0){13}}
\put(114,12){\line(0,1){7}}
\put(101,19){\line(1,0){13}}
\put(125,19){102}
\put(121.5,22){$\langle3,K_4,9\rangle$}
\put(121,18){\line(0,1){7}}
\put(121,18){\line(1,0){13}}
\put(134,18){\line(0,1){7}}
\put(121,25){\line(1,0){13}}
\put(15,25){\vector(3,1){45}}
\put(15,25){\vector(1,0){25}}
\put(15,25){\vector(1,-4){5}}
\put(18,28){$\pi_3$}
\put(38,7){$\pi_1$}
\put(35,4){\vector(1,0){45}}
\put(35,4){\line(5,1){50}}
\put(85,14){\vector(1,0){15}}
\put(35,4){\line(3,1){54}}
\put(89,22){\vector(1,0){31}}
\put(120,0){\footnotesize{\textbf{(c)}}}
\end{picture}
\end{center}
\normalsize
\caption{(a)-(c) are square trees from roots $\langle1,K_6,1\rangle$, $\langle2,K_7,1\rangle$, $\langle3,K_8,1\rangle$, respectively.}
\label{fig:1}
\end{figure}
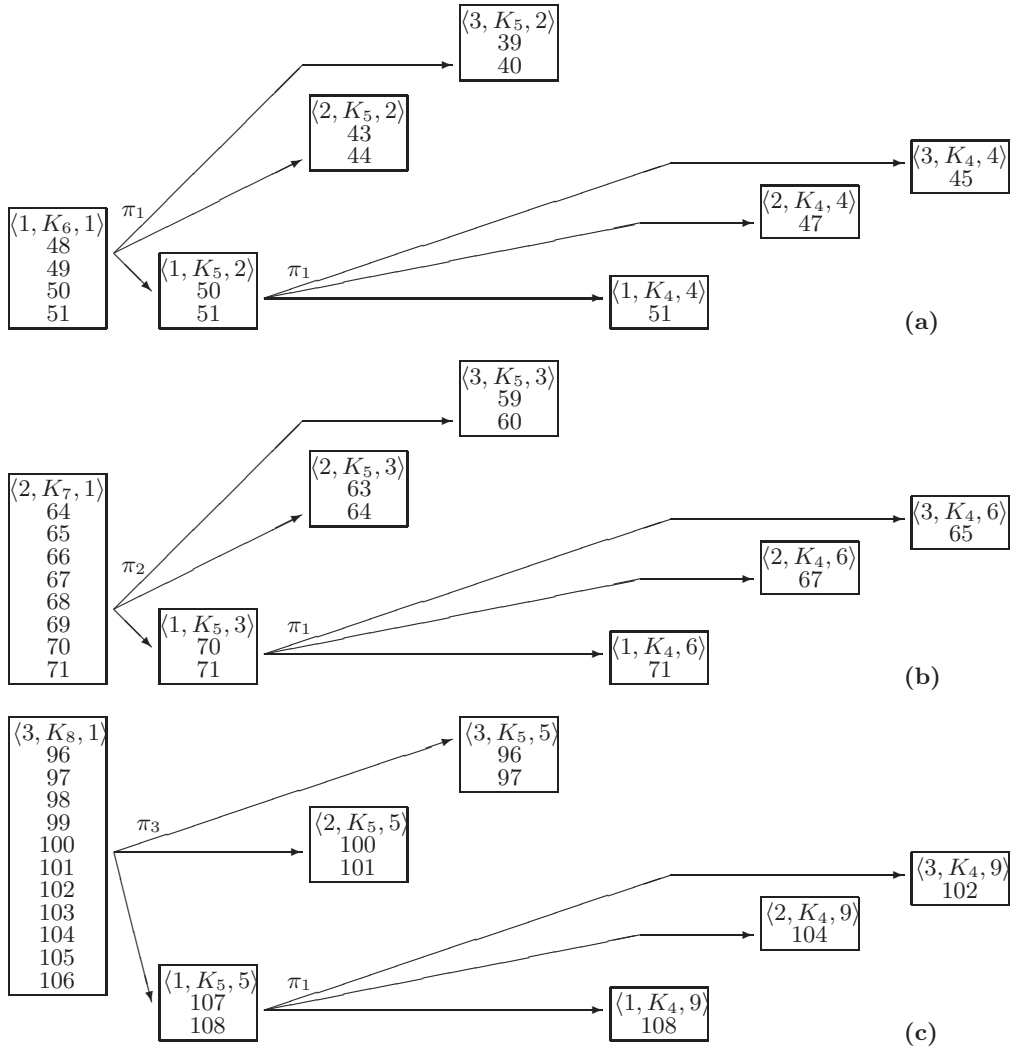

On the other hand,
for $m\geq4$ and $j\in\{1,2,3\}$, each $\langle j,K_{m},1\rangle$ belongs to the square trees. Moreover $\langle j,K_m,\hat{a}\rangle$ (resp. $\langle j,K_m,\hat{b}\rangle$, $\langle j,K_m,\hat{c}\rangle$) is subset of $\pi_1\langle1,K_{m+1},p\rangle$
(resp. $\pi_2\langle2,K_{m+2},p\rangle$, $\pi_3\langle3,K_{m+3},p\rangle$).
Notice that $\mathbb{N}=\{1\}\cup\{P(a,p)+1\}\cup\{P(b,p)+1\}\cup\{P(c,p)+1\}$, the square trees contain all $\langle j,K_m,p\rangle$, i.e. the positions of all squares in $\mathbb{T}$.
Figure \ref{fig:1} shows some examples.

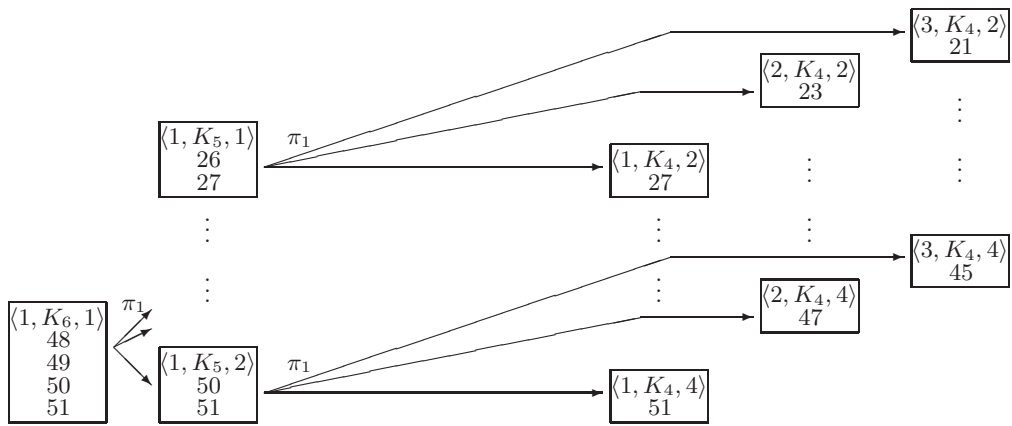
\begin{figure}
\footnotesize
\setlength{\unitlength}{1mm}
\begin{center}
\begin{picture}(140,54)
\put(6,1){51}
\put(6,4){50}
\put(6,7){49}
\put(6,10){48}
\put(1,13){$\langle1,K_6,1\rangle$}
\put(1,0){\line(0,1){16}}
\put(1,0){\line(1,0){13}}
\put(14,0){\line(0,1){16}}
\put(1,16){\line(1,0){13}}
\put(15,10){\vector(1,-1){5}}
\put(15,10){\vector(2,1){5}}
\put(15,10){\vector(1,1){5}}
\put(16,15){$\pi_1$}
\put(26,1){51}
\put(26,4){50}
\put(21,7){$\langle1,K_5,2\rangle$}
\put(21,0){\line(0,1){10}}
\put(21,0){\line(1,0){13}}
\put(34,0){\line(0,1){10}}
\put(21,10){\line(1,0){13}}
\put(86,1){51}
\put(81,4){$\langle1,K_4,4\rangle$}
\put(81,0){\line(0,1){7}}
\put(81,0){\line(1,0){13}}
\put(94,0){\line(0,1){7}}
\put(81,7){\line(1,0){13}}
\put(106,13){47}
\put(101,16){$\langle2,K_4,4\rangle$}
\put(101,12){\line(0,1){7}}
\put(101,12){\line(1,0){13}}
\put(114,12){\line(0,1){7}}
\put(101,19){\line(1,0){13}}
\put(126,19){45}
\put(121,22){$\langle3,K_4,4\rangle$}
\put(121,18){\line(0,1){7}}
\put(121,18){\line(1,0){13}}
\put(134,18){\line(0,1){7}}
\put(121,25){\line(1,0){13}}
\put(38,7){$\pi_1$}
\put(35,4){\vector(1,0){45}}
\put(35,4){\line(5,1){50}}
\put(85,14){\vector(1,0){15}}
\put(35,4){\line(3,1){54}}
\put(89,22){\vector(1,0){31}}
\put(26,31){27}
\put(26,34){26}
\put(21,37){$\langle1,K_5,1\rangle$}
\put(21,30){\line(0,1){10}}
\put(21,30){\line(1,0){13}}
\put(34,30){\line(0,1){10}}
\put(21,40){\line(1,0){13}}
\put(86,31){27}
\put(81,34){$\langle1,K_4,2\rangle$}
\put(81,30){\line(0,1){7}}
\put(81,30){\line(1,0){13}}
\put(94,30){\line(0,1){7}}
\put(81,37){\line(1,0){13}}
\put(106,43){23}
\put(101,46){$\langle2,K_4,2\rangle$}
\put(101,42){\line(0,1){7}}
\put(101,42){\line(1,0){13}}
\put(114,42){\line(0,1){7}}
\put(101,49){\line(1,0){13}}
\put(126,49){21}
\put(121,52){$\langle3,K_4,2\rangle$}
\put(121,48){\line(0,1){7}}
\put(121,48){\line(1,0){13}}
\put(134,48){\line(0,1){7}}
\put(121,55){\line(1,0){13}}
\put(38,37){$\pi_1$}
\put(35,34){\vector(1,0){45}}
\put(35,34){\line(5,1){50}}
\put(85,44){\vector(1,0){15}}
\put(35,34){\line(3,1){54}}
\put(89,52){\vector(1,0){31}}
\put(27,24){$\vdots$}
\put(27,16){$\vdots$}
\put(87,24){$\vdots$}
\put(87,16){$\vdots$}
\put(107,24){$\vdots$}
\put(107,32){$\vdots$}
\put(127,32){$\vdots$}
\put(127,40){$\vdots$}
\end{picture}
\end{center}
\normalsize
\caption{An example of the graph embedding in square trees.}
\label{fig:2}
\end{figure}

By the square trees, we have the relation between the number of squares ending at position $\Gamma_{j,m,p}[i]$ and $\Gamma_{j,m,1}[i]$, see Property \ref{P5.5}.
Figure \ref{fig:2} shows the relation: in the tree with root $\langle1,K_6,1\rangle$,
the branch from node $\langle1,K_5,2\rangle$ is the graph embedding of the tree with root $\langle1,K_5,1\rangle$.

\begin{property}[]\label{P5.5} For $j=1,2,3$, $m\geq4$ and $p\geq1$,
$$\{\omega:\omega\omega\triangleright\mathbb{T}[1,\Gamma_{j,m,1}[i]]\}
=\{\omega:\omega\omega\triangleright\mathbb{T}[1,\Gamma_{j,m,p}[i]],Ker(\omega)=K_h,4\leq h\leq m\}.$$
Here $1\leq i\leq k_{m}-1$ for $j=1,2$; $1\leq i\leq t_{m-4}-k_{m-3}+1$ for $j=3$.
\end{property}

\noindent\emph{Example.} Let $j=1$, $m=5$, $p=3$, $i=2$. All squares ending at position $\Gamma_{1,5,1}[2]=27$ are $\{\omega\omega,\varpi\varpi\}$, where $\omega=abacaba$ and $\varpi=bacabaabacaba$.
All squares ending at position
$\Gamma_{1,5,3}[2]=71$ are $\{\omega\omega,\varpi\varpi,\mu\mu\}$, where $\mu=abacababacabaabacaba$.
Since $Ker(\omega\omega)=aa=K_{4}$, $Ker(\varpi\varpi)=bab=K_5$ and $Ker(\mu\mu)=aabacabaa=K_7$,
only $\{\omega\omega,\varpi\varpi\}$ are squares with kernel $K_{h}$, $4\leq h\leq 5$.

\vspace{0.2cm}

Taking $p=1$, we have
\begin{equation*}
\begin{cases}
\Gamma_{1,m,1}=[\frac{t_{m}+2t_{m-1}-t_{m-2}-1}{2},\cdots,
\frac{t_{m}+2t_{m-1}+t_{m-2}-3}{2}];\\
\Gamma_{2,m,1}=[\frac{-t_{m}+4t_{m-1}+t_{m-2}-1}{2},\cdots,
\frac{t_{m}+2t_{m-1}-t_{m-2}-3}{2}];\\
\Gamma_{3,m,1}=[\frac{t_{m}+t_{m-2}-1}{2},\cdots,
\frac{-t_{m}+4t_{m-1}+t_{m-2}-3}{2}].
\end{cases}
\end{equation*}

Denote $b(n)=\sharp\{(\omega,p):\omega_p\omega_{p+1}\triangleright\mathbb{T}[1,n]\}$
the number of squares ending at position $n$.
By Property \ref{b} below, we can calculate $b(n)$, and obversely calculate $B(n)$ by $B(n)=\sum_{i=1}^n b(i)$.

\begin{property}\label{b}\ $b(i)=0$ for $i\leq7$, $b([8])=[1]$, $b([9,10])=[0,1]$, $b([11,\cdots,14])=[0,0,0,1]$,
$b([15,16])=[1,1]$, $b([17,\cdots,20])=[0,0,1,1]$, $b([28,\cdots,31])=[1,1,1,1]$,
\begin{equation*}
\begin{cases}
b(\Gamma_{1,m,1})
=[b(\Gamma_{3,m-1,1}),b(\Gamma_{2,m-1,1}),b(\Gamma_{1,m-1,1})]
+[\underbrace{0,\cdots,0}_{t_{m-2}-k_{m}+1},\underbrace{1,\cdots,1}_{k_{m}-1}]
\text{ for }m\geq5;\\
b(\Gamma_{2,m,1})
=[b(\Gamma_{3,m-2,1}),b(\Gamma_{2,m-2,1}),b(\Gamma_{1,m-2,1})]
+[\underbrace{0,\cdots,0}_{t_{m-3}-k_{m}+1},\underbrace{1,\cdots,1}_{k_{m}-1}]
\text{ for }m\geq6;\\
b(\Gamma_{3,m,1})
=[b(\Gamma_{3,m-3,1}),b(\Gamma_{2,m-3,1}),b(\Gamma_{1,m-3,1})]
+[\underbrace{1,\cdots,1}_{t_{m-4}-k_{m-3}+1},\underbrace{0,\cdots,0}_{k_{m-3}-1}]
\text{ for }m\geq7.
\end{cases}
\end{equation*}
\end{property}

Obversely we can calculate $B(n)$ by $B(n)=\sum_{i=1}^n b(i)$. But when $n$ large, the method is complicated. Now we turn to give a fast algorithm.
Denote the sum of all elements in vector $b(\Gamma_{j,m,1})$ by $\sum b(\Gamma_{j,m,1})$ for $j\in\{1,2,3\}$.
Denote $\Phi_m=\sum b(\Gamma_{3,m,1})+\sum b(\Gamma_{2,m,1})+\sum b(\Gamma_{1,m,1})$.
The immediately corollaries of Property \ref{b} are
$\sum b(\Gamma_{1,m,1})=\Phi_{m-1}+k_m-1$,
$\sum b(\Gamma_{2,m,1})=\Phi_{m-2}+k_m-1$,
$\sum b(\Gamma_{3,m,1})=\Phi_{m-3}+t_{m-4}-k_{m-3}+1$.
Moreover,
$\Phi_m=\Phi_{m-1}+\Phi_{m-2}+\Phi_{m-3}+\tfrac{-3t_{m}+6t_{m-1}+t_{m-2}-1}{2}$ for $m\geq7$.
By induction, we can prove Property \ref{P4.2} below easily.

\begin{property}[]\label{P4.2}
For $m\geq4$, (1) $\Phi_m=\tfrac{m}{22}(-5t_{m}+14t_{m-1}+4t_{m-2})
+\tfrac{1}{44}(67t_{m}-166t_{m-1}+5t_{m-2})+\tfrac{1}{4}$;

(2) $\sum b(\Gamma_{1,m,1})=\frac{m}{22}(4t_{m}-9t_{m-1}+10t_{m-2})
+\frac{1}{44}(19t_{m}+36t_{m-1}-169t_{m-2})-\frac{1}{4}$,

$\sum b(\Gamma_{2,m,1})=\frac{m}{22}(10t_{m}-6t_{m-1}-19t_{m-2})
+\frac{1}{44}(-189t_{m}+156t_{m-1}+331t_{m-2})-\frac{1}{4}$,

$\sum b(\Gamma_{3,m,1})=\frac{m}{22}(-19t_{m}+29t_{m-1}+13t_{m-2})
+\frac{1}{44}(237t_{m}-358t_{m-1}-157t_{m-2})+\frac{3}{4}$;

(3) $\sum_{j=4}^{m-1}\Phi_j=\tfrac{m}{44}(13t_{m}-10t_{m-1}+5t_{m-2})
+\tfrac{2}{11}(-8t_{m}+8t_{m-1}-7t_{m-2})+\tfrac{m}{4}+2$;

(4) $B(\max\Gamma_{3,m,1})=\frac{m}{44}(-25t_{m}+48t_{m-1}+31t_{m-2})
+\frac{1}{44}(173t_{m}-294t_{m-1}-213t_{m-2})+\frac{m+11}{4}$,

$B(\max\Gamma_{2,m,1})=\frac{m}{44}(-5t_{m}+36t_{m-1}-7t_{m-2})
+\frac{1}{22}(-8t_{m}-69t_{m-1}+59t_{m-2})+\frac{m+10}{4}$,

$B(\max\Gamma_{1,m,1})=\frac{m}{44}(3t_{m}+18t_{m-1}+13t_{m-2})
+\frac{1}{44}(3t_{m}-102t_{m-1}-51t_{m-2})+\frac{m+9}{4}$.
\end{property}

For any $n\geq52$, let $m$ such that
$n\in[\Gamma_{3,m,1},\Gamma_{2,m,1},\Gamma_{1,m,1}]
=[\tfrac{t_{m}+t_{m-2}-1}{2},\cdots,\tfrac{t_{m}+2t_{m-1}+t_{m-2}-3}{2}].$
We already determine the expression of $B(\max{\Gamma_{j,m,1}})$ for $j\in\{1,2,3\}$, $m\geq5$. In order to calculate $B(n)$, we only need to calculate $\sum_{i=\min{\Gamma_{j,m,1}}}^n b(i)$.

\textbf{Case 1.} When $n\in\Gamma_{3,m,1}=[\frac{t_{m}+t_{m-2}-1}{2},\cdots,
\frac{-t_{m}+4t_{m-1}+t_{m-2}-3}{2}]$ for $m\geq7$. Denote
\begin{equation*}
\begin{cases}
\theta_m^1=\min\Gamma_{3,m,1}=\frac{t_{m}+t_{m-2}-1}{2};\\
\theta_m^2=\min\Gamma_{3,m,1}+|\Gamma_{3,m-3,1}|=\frac{-5t_{m}+10t_{m-1}+3t_{m-2}-1}{2};\\
\theta_m^3=\min\Gamma_{3,m,1}+|\Gamma_{3,m-3,1}|+|\Gamma_{2,m-3,1}|
=\frac{-t_{m}+6t_{m-1}-3t_{m-2}-1}{2};\\
\eta_m^1=\min\Gamma_{3,m,1}+t_{m-4}-k_{m-3}+1=-2t_{m}+5t_{m-1}.\\
\theta_m^4=\max\Gamma_{3,m,1}+1=\min\Gamma_{2,m,1}=\frac{-t_{m}+4t_{m-1}+t_{m-2}-1}{2}.
\end{cases}
\end{equation*}
Obviously, $\theta_m^3<\eta_m^1<\theta_m^4$ for $m\geq7$,
$\min\Gamma_{3,m,1}-\min\Gamma_{3,m-3,1}=t_{m-1}$. By Property \ref{b}, we have:

\begin{property}[]\label{P4.3}\
For $n\geq52$, let $m$ such that $n\in\Gamma_{3,m,1}$, then $m\geq7$ and
\begin{equation*}
\begin{split}
&\sum\nolimits_{i=\min\Gamma_{3,m,1}}^nb(i)\\
=&\begin{cases}
\sum_{i=\min\Gamma_{3,m-3,1}}^{n-t_{m-1}}b(i)+n-\min\Gamma_{3,m,1}+1,
&\theta_m^1\leq n<\theta_m^2;\\
\sum_{i=\min\Gamma_{2,m-3,1}}^{n-t_{m-1}}b(i)+\sum b(\Gamma_{3,m-3,1})
+n-\min\Gamma_{3,m,1}+1,&\theta_m^2\leq n<\theta_m^3;\\
\sum_{i=\min\Gamma_{1,m-3,1}}^{n-t_{m-1}}b(i)+\sum b(\Gamma_{3,m-3,1})
+\sum b(\Gamma_{2,m-3,1})+n-\min\Gamma_{3,m,1}+1,&\theta_m^3\leq n<\eta_m^1;\\
\sum_{i=\min\Gamma_{1,m-3,1}}^{n-t_{m-1}}b(i)+\sum b(\Gamma_{3,m-3,1})
+\sum b(\Gamma_{2,m-3,1})+\frac{-5t_{m}+10t_{m-1}-t_{m-2}+1}{2},&otherwise.
\end{cases}
\end{split}
\end{equation*}
\end{property}

\textbf{Case 2.} When $n\in\Gamma_{2,m,1}=[\frac{-t_{m}+4t_{m-1}+t_{m-2}-1}{2},\cdots,
\frac{t_{m}+2t_{m-1}-t_{m-2}-3}{2}]$ for $m\geq6$. Denote
\begin{equation*}
\begin{cases}
\theta_m^5=\min\Gamma_{2,m,1}+|\Gamma_{3,m-2,1}|=\frac{3t_{m}-5t_{m-2}-1}{2};\\
\eta_m^2=\min\Gamma_{2,m,1}+t_{m-3}-k_{m}+1=2t_{m-1}-t_{m-2}.\\
\theta_m^6=\min\Gamma_{3,m,1}+|\Gamma_{3,m-3,1}|+|\Gamma_{2,m-3,1}|
=\frac{3t_{m}-2t_{m-1}-t_{m-2}-1}{2};\\
\theta_m^7=\max\Gamma_{2,m,1}+1=\min\Gamma_{1,m,1}=\frac{t_{m}+2t_{m-1}-t_{m-2}-1}{2}.
\end{cases}
\end{equation*}
Obviously, $\theta_m^5<\eta_m^2\leq\theta_m^6$ for $m\geq6$, $\min\Gamma_{2,m,1}-\min\Gamma_{3,m-2,1}=t_{m-1}$. By Property \ref{b}, we have:

\begin{property}[]\label{P4.4}\
For $n\geq32$, let $m$ such that $n\in\Gamma_{2,m,1}$, then $m\geq6$ and
\begin{equation*}
\begin{split}
&\sum\nolimits_{i=\min\Gamma_{2,m,1}}^nb(i)\\
=&\begin{cases}
\sum_{i=\min\Gamma_{3,m-2,1}}^{n-t_{m-1}}b(i),
&\theta_m^4\leq n<\theta_m^5;\\
\sum_{i=\min\Gamma_{2,m-2,1}}^{n-t_{m-1}}b(i)+\sum b(\Gamma_{3,m-2,1}),
&\theta_m^5\leq n<\eta_m^2;\\
\sum_{i=\min\Gamma_{2,m-2,1}}^{n-t_{m-1}}b(i)+\sum b(\Gamma_{3,m-2,1})
+n-\eta_m^2+1,&\eta_m^2\leq n<\theta_m^6;\\
\sum_{i=\min\Gamma_{1,m-2,1}}^{n-t_{m-1}}b(i)+\sum b(\Gamma_{3,m-2,1})
+\sum b(\Gamma_{2,m-2,1})+n-\eta_m^2+1,&otherwise.
\end{cases}
\end{split}
\end{equation*}
\end{property}

\textbf{Case 3.} When $n\in\Gamma_{1,m,1}=[\frac{t_{m}+2t_{m-1}-t_{m-2}-1}{2},\cdots,
\frac{t_{m}+2t_{m-1}+t_{m-2}-3}{2}]$ for $m\geq5$. Denote
\begin{equation*}
\begin{cases}
\theta_m^8=\min\Gamma_{1,m,1}+|\Gamma_{3,m-1,1}|=\frac{t_{m}+3t_{m-2}-1}{2};\\
\theta_m^9=\min\Gamma_{1,m,1}+|\Gamma_{3,m-1,1}|+|\Gamma_{2,m-1,1}|
=\frac{-t_{m}+4t_{m-1}+3t_{m-2}-1}{2};\\
\eta_m^3=\min\Gamma_{1,m,1}+t_{m-2}-k_{m}+1=2t_{m-1}.\\
\theta_m^{10}=\max\Gamma_{1,m,1}+1=\min\Gamma_{3,m+1,1}=\frac{t_{m}+2t_{m-1}+t_{m-2}-1}{2}.
\end{cases}
\end{equation*}
Obviously, $\theta_m^9<\eta_m^3<\theta_m^{10}$ for $m\geq5$, $\min\Gamma_{1,m,1}-\min\Gamma_{3,m-1,1}=t_{m-1}$. By Property \ref{b}, we have:

\begin{property}[]\label{P4.5}\
For $n\geq21$, let $m$ such that $n\in\Gamma_{1,m,1}$, then $m\geq5$ and
\begin{equation*}
\begin{split}
&\sum\nolimits_{i=\min\Gamma_{1,m,1}}^nb(i)\\
=&\begin{cases}
\sum_{i=\min\Gamma_{3,m-1,1}}^{n-t_{m-1}}b(i),
&\theta_m^7\leq n<\theta_m^8;\\
\sum_{i=\min\Gamma_{2,m-1,1}}^{n-t_{m-1}}b(i)+\sum b(\Gamma_{3,m-1,1}),
&\theta_m^8\leq n<\theta_m^9;\\
\sum_{i=\min\Gamma_{1,m-1,1}}^{n-t_{m-1}}b(i)+\sum b(\Gamma_{3,m-1,1})
+\sum b(\Gamma_{2,m-1,1}),&\theta_m^9\leq n<\eta_m^3;\\
\sum_{i=\min\Gamma_{1,m-1,1}}^{n-t_{m-1}}b(i)+\sum b(\Gamma_{3,m-1,1})
+\sum b(\Gamma_{2,m-1,1})+n-\eta_m^3+1,&otherwise.
\end{cases}
\end{split}
\end{equation*}
\end{property}

\begin{aB}[]
Step 1. For $n\leq51$, calculate $\sum_{i=1}^n b(i)$ by Property \ref{b}.

Step 2. For $n\geq52$, find the $m$ and $j$ such that $n\in\Gamma_{j,m,1}$, then $m\geq7$.
Calculate $B(\min{\Gamma_{j,m,1}}-1)$ by Property \ref{P4.2}.
And calculate $\sum_{i=\min{\Gamma_{j,m,1}}}^n b(i)$ by Property \ref{P4.3}-\ref{P4.5}.

Step 3. $B(n)=B(\min{\Gamma_{j,m,1}}-1)+\sum_{i=\min{\Gamma_{j,m,1}}}^n b(i)$.
\end{aB}

\noindent\emph{Example.} Consider $n=60\in\Gamma_{2,7,1}=[59,\cdots,71]$.
$\theta_7^4=\frac{-t_{7}+4t_{6}+t_{5}-1}{2}=59\leq n
<\theta_7^5=\frac{3t_{7}-5t_{5}-1}{2}=61$.
By Property \ref{P4.4},
$\sum_{i=\min\Gamma_{2,7,1}}^{60}b(i)
=\sum_{i=\min\Gamma_{3,5,1}}^{60-t_{6}}b(i)
=\sum_{i=15}^{16}b(i)=b(15)+b(16)=2$.
By Property \ref{P4.2},
$B(\min{\Gamma_{2,7,1}}-1)=B(\max{\Gamma_{3,7,1}})
=\frac{7}{44}(-25t_{7}+48t_{6}+31t_{5})
+\frac{1}{44}(173t_{7}-294t_{6}-213t_{5})+\frac{18}{4}=45$.
Thus $B(60)=B(\max{\Gamma_{3,7,1}})+\sum_{i=15}^{16} b(i)=47$.

\vspace{0.2cm}

Now we turn to give the expressions of $B(t_m)$.

For $m\geq7$, $\theta_m^8\leq t_{m}<\theta_m^9$
and $\theta_{m-1}^6\leq t_{m}-t_{m-1}<\theta_{m-1}^7$. By Property \ref{P4.4}-\ref{P4.5},
$$\begin{array}{rl} &\sum_{i=\min\Gamma_{1,m,1}}^{t_m}b(i)-\sum_{i=\min\Gamma_{1,m-3,1}}^{t_{m-3}}b(i)\\
=&\sum b(\Gamma_{3,m-1,1})+\sum b(\Gamma_{3,m-3,1})+\sum b(\Gamma_{2,m-3,1})
+2t_{m}-2t_{m-1}-3t_{m-2}+1\\
=&\tfrac{m}{22}(-19t_m+29t_{m-1}+13t_{m-2})
+\tfrac{1}{44}(347t_m-622t_{m-1}-47t_{m-2})+\tfrac{9}{4}.
\end{array}$$
For $m\geq4$, by induction, $\sum_{i=\min\Gamma_{1,m,1}}^{t_m}b(i)$ is equal to
$$\tfrac{m}{44}(23t_m-38t_{m-1}-3t_{m-2})
+\tfrac{1}{44}(-65t_m+164t_{m-1}-105t_{m-2})+\tfrac{3m}{4}-\tfrac{9}{4}.$$
Since $\min\Gamma_{1,m,1}-1=\max\Gamma_{2,m,1}$, $B(t_m)=B(\max\Gamma_{2,m,1})+\sum_{i=\min\Gamma_{1,m,1}}^{t_m}b(i)$.
By Property \ref{P4.2}, we can prove Theorem 21 in H.Mousavi and J.Shallit\cite{MS2014} in a novel way, i.e. \textbf{Theorem B}.

\section{The number of distinct cubes, $C(n)$}

Let $\omega$ be a factor with kernel $K_m$. By an analogous argument as Section 3,
$\omega_{p}\omega_{p+1}\omega_{p+2}\prec\mathbb{T}$ has only one case:
$G_p(K_m)=G_{p+1}(K_m)=G_1(K_m)$. Here $|G_1(K_m)|=t_{m}-k_{m}$, and
$$\omega\omega\omega
=T_{m-1}[i,t_{m-1}]T_{m-2}T_{m-3}[1,t_{m-3}-1]\underline{K_{m+4}}T_{m+1}[k_{m+4},t_{m+1}]T_{m-3}T_{m-2}[1,i-1],$$
where $1\leq i\leq k_{m+1}-1$, $m\geq3$, $Ker(\omega\omega\omega)=K_{m+4}$ and $|\omega|=t_{m}$.

\vspace{0.2cm}

\noindent\emph{Remark.}
By this case, we have that: all cubes in $\mathbb{T}$ are of length $3t_m$ for some $m\geq3$. Furthermore, for all $m\geq3$, there exists a cube of length $3t_m$ in $\mathbb{T}$. This is Theorem 7 in \cite{MS2014}.

\vspace{0.2cm}

By the case of cubes, we define a set for $m\geq7$,
\begin{equation*}
\begin{split}
\langle K_m,p\rangle
=&\{P(\omega\omega\omega,p):Ker(\omega\omega\omega)=K_m,|\omega|=t_{m-4},\omega\omega\omega\prec\mathbb{T}\}\\
=&\{P(K_m,p)+\tfrac{-t_{m-2}+5t_{m-4}+1}{2}
,\cdots,P(K_m,p)+t_{m-4}-1\}.
\end{split}
\end{equation*}
Then $\sharp\langle K_m,p\rangle=k_{m-3}-1$.
$P(K_m,1)=\frac{t_{m}+t_{m-2}-1}{2}$, so
$\langle K_m,1\rangle=\{t_{m-1}+2t_{m-4}
,\cdots,\frac{3t_{m-1}-t_{m-3}-3}{2}\}$.
Since $\langle K_m,1\rangle$ are pairwise disjoint for different $m$, we get a chain
$\langle K_7,1\rangle,\langle K_8,1\rangle,\cdots,
\langle K_m,1\rangle,\cdots$.
Denote $c(n)=\sharp\{\omega\omega\omega:\omega\omega\omega\triangleright\mathbb{T}[1,n],
\omega\omega\omega\not\!\prec\mathbb{T}[1,n-1]\}$,
then $c(n)=1$ if and only if $n\in\cup_{m\geq7}\langle K_m,1\rangle$.

\begin{property}[]\label{c} $c(n)=0$ for $n\leq57$; for $n\geq58$, let $m$ such that $t_{m-1}+2t_{m-4}\leq n<t_{m}+2t_{m-3}$,
then $m\geq7$ and $c(n)=1$ if and only if $n\leq t_{m-1}+k_{m+1}-2$.
\end{property}

Obviously, $C(n)=\sum_{i=1}^n c(i)$, and $C(n)=0$ for $n\leq57$.
For $m\geq7$, we denote
$$\hat{\alpha}_m=\min\langle K_m,1\rangle=t_{m-1}+2t_{m-4},~
\hat{\beta}_m=\max\langle K_m,1\rangle
=\tfrac{3t_{m-1}-t_{m-3}-3}{2}.$$
When $n\geq58$, find $m$ such that $\hat{\alpha}_m\leq n<\hat{\alpha}_{m+1}$.
By Property \ref{c} and the definition of $\Delta_m$,
\begin{equation*}
\begin{cases}
C(\hat{\alpha}_m)=C(\min\langle K_m,1\rangle)=\Delta_{m-4}+1=\frac{t_{m-6}+t_{m-7}-m+6}{2},\\
C(\hat{\beta}_m)=C(\max\langle K_{m},1\rangle)=\Delta_{m-3}=\frac{t_{m-5}+t_{m-6}-m+3}{2}.
\end{cases}
\end{equation*}
Moreover, (1) when $\hat{\alpha}_m\leq n\leq\hat{\beta}_m$,
$C(n)=C(\hat{\alpha}_m)+n-\hat{\alpha}_m$;
(2) when $\hat{\beta}_m< n<\hat{\alpha}_{m+1}$,
$C(n)=C(\hat{\beta}_m)$.
Thus we get \textbf{Theorem C.}, i.e. the expression of $C(n)$ for all $n$.

\vspace{0.2cm}

For $m\geq7$, since $t_m\geq\hat{\beta}_{m}=t_{m-1}+k_{m+1}-2$, by Theorem C., we have
$C(t_m)=C(\hat{\beta}_m)$, i.e.

\begin{theorem}[]\ For $m\leq6$, $C(t_m)=0$, for $m\geq7$,
$C(t_m)=\frac{t_{m-5}+t_{m-6}-m+3}{2}$.
\end{theorem}

\section{The number of repeated cubes, $D(n)$}

For $m\geq7$ and $p\geq1$, we consider the vectors $\Gamma_{m,p}=[P(K_m,p),\cdots,P(K_{m},p)+t_{m-1}-1]$.
Using Property \ref{L},
comparing minimal and maximal elements in these sets below, we have
$$\Gamma_{m,p}=[\Gamma_{m-3,P(a,p)+1},\Gamma_{m-2,P(b,p)+1},\Gamma_{m-1,P(c,p)+1}]
\text{ for }m\geq10.$$
Thus we establish the recursive relations for all $\Gamma_{m,p}$, $m\geq10$.
We arrange all elements of $\langle K_m,p\rangle$ in ascending order as vector $[\langle K_m,p\rangle]$. It is easy to check that
$$\Gamma_{m,p}
=[\underbrace{P,\cdots,P+\tfrac{-t_{m-2}+5t_{m-4}-1}{2}}_{\tfrac{-t_{m-2}+5t_{m-4}+1}{2}},\underbrace{[\langle K_m,p\rangle]}_{\tfrac{t_{m-2}-3t_{m-4}-1}{2}},
\underbrace{P+t_{m-4},\cdots,P+t_{m-1}-1}_{t_{m-2}+t_{m-3}}].$$
Here we denote $P(K_m,p)$ by $P$ for short. The number under ``$\stackrel{\underbrace{}}{}$" means the number of elements.

By an analogous argument in Section 4,
we get the recursive structure of the positions of repeated cubes in $\mathbb{T}$, called cube trees. Here $\pi$ is a substitution over $\{\langle K_m,p\rangle:m\geq7,p\geq1\}$ that
\begin{equation*}
\begin{cases}
\pi\langle K_8,p\rangle=\langle K_7,P(a,p)+1\rangle;~
\pi\langle K_9,p\rangle=\langle K_7,P(b,p)+1\rangle\cup\langle K_8,P(a,p)+1\rangle;\\
\pi\langle K_m,p\rangle=\langle K_{m-3},P(a,p)+1\rangle\cup\langle K_{m-2},P(b,p)+1\rangle
\cup\langle K_{m-1},P(c,p)+1\rangle,m\geq10.
\end{cases}
\end{equation*}
The cube trees contain all $\langle K_m,p\rangle$, i.e. the positions of all cubes in $\mathbb{T}$.

Let $p=1$, $\Gamma_{m,1}=[P(K_m,1),\cdots,P(K_{m+1},1)-1]
=[\frac{t_{m}+t_{m-2}-1}{2},\cdots,\frac{t_{m+1}+t_{m-1}-3}{2}]$.

Denote $d(n)=\sharp\{(\omega,p):\omega_p\omega_{p+1}\omega_{p+2}\triangleright\mathbb{T}[1,n]\}$
the number of cubes ending at position $n$.
By Property \ref{d} below, we can calculate $d(n)$, and obversely calculate $D(n)$ by $D(n)=\sum_{i=1}^n d(i)$.

\begin{property}\label{d}\ For $n\leq51$, $d(n)=0$; and
\begin{equation*}
\begin{cases}
d(\Gamma_{7,1})=d([52,\cdots,95])=[\underbrace{0,\cdots,0}_6,1,\underbrace{0,\cdots,0}_{37}];\\
d(\Gamma_{8,1})=d([96,\cdots,176])
=[\underbrace{0,\cdots,0}_{11},1,1,\underbrace{0,\cdots,0}_{30},1,\underbrace{0,\cdots,0}_{37}];\\
d(\Gamma_{9,1})=d([177,\cdots,325])
=[\underbrace{0,\cdots,0}_{20},\underbrace{1,\cdots,1}_4,\underbrace{0,\cdots,0}_{6},1,
\underbrace{0,\cdots,0}_{48},1,1,\underbrace{0,\cdots,0}_{30},1,\underbrace{0,\cdots,0}_{37}];\\
d(\Gamma_{m,1})=[d(\Gamma_{m-3,1}),d(\Gamma_{m-2,1}),d(\Gamma_{m-1,1})]+
[\underbrace{0,\cdots,0,}_{\tfrac{-t_{m-2}+5t_{m-4}+1}{2}}
\underbrace{1,\cdots,1,}_{\tfrac{t_{m-2}-3t_{m-4}-1}{2}}
\underbrace{0,\cdots,0}_{t_{m-2}+t_{m-3}}], m\geq10.
\end{cases}
\end{equation*}
\end{property}

Now we turn to give a fast algorithm.
Denote the sum of all elements in vector $d(\Gamma_{m,1})$ by $\sum d(\Gamma_{m,1})$.
An immediately corollaries of Property \ref{d} are
$\sum d(\Gamma_{m,1})=\sum d(\Gamma_{m-3,1})+\sum d(\Gamma_{m-2,1})+\sum d(\Gamma_{m-1,1})+\tfrac{t_{m-2}-3t_{m-4}-1}{2}$.
By induction, we can prove Property \ref{P6.2} below easily.

\begin{property}[]\label{P6.2} For $m\geq7$,
(1) $D(\min\Gamma_{m+1,1})=D(\max\Gamma_{m,1})$;

(2) $\sum d(\Gamma_{m,1})=\frac{m}{22}(7t_{m}-13t_{m-1}+t_{m-2})
+\frac{1}{44}(-41t_{m}+74t_{m-1}-7t_{m-2})+\frac{1}{4}$;

(3) $D(\max\Gamma_{m,1})=\sum_{j=7}^{m}\sum d(\Gamma_{j,1})
=\frac{m}{44}(9t_{m}-12t_{m-1}-5t_{m-2})
+\frac{3}{11}(-2t_{m}+2t_{m-1}+t_{m-2})+\frac{m}{4}$.
\end{property}

For any $n\geq52$, let $m$ such that
$n\in\Gamma_{m,1}$. Now we calculate $\sum_{i=\min{\Gamma_{m,1}}}^n d(i)$.
Denote
\begin{equation*}
\begin{cases}
\hat{\theta}_m^1=\min\Gamma_{m,1}=\tfrac{t_{m}+t_{m-2}-1}{2};\\
\hat{\eta}_m^1=\min\Gamma_{m,1}+\tfrac{-t_{m-2}+5t_{m-4}+1}{2}=\tfrac{t_{m}+5t_{m-4}}{2};\\
\hat{\eta}_m^2=\min\Gamma_{m,1}+\tfrac{-t_{m-2}+5t_{m-4}+1}{2}+\tfrac{t_{m-2}-3t_{m-4}-1}{2}=\tfrac{3t_{m-1}-t_{m-3}-1}{2};\\
\hat{\theta}_m^2=\min\Gamma_{m,1}+|\Gamma_{m-3,1}|=\tfrac{3t_{m-1}-t_{m-3}-1}{2}=\hat{\eta}_m^2;\\
\hat{\theta}_m^3=\min\Gamma_{m,1}+|\Gamma_{m-3,1}|+|\Gamma_{m-2,1}|=\tfrac{3t_{m-1}+t_{m-3}-1}{2}.
\end{cases}
\end{equation*}
Obviously, $\hat{\theta}_m^1<\hat{\eta}_m^1<\hat{\eta}_m^2=\hat{\theta}_m^2<\hat{\theta}_m^3<\hat{\theta}_{m+1}^1$ for $m\geq7$. By Property \ref{d}, we have:

\begin{property}[]\label{P6.3}\
For $n\geq326$, let $m$ such that $n\in\Gamma_{m,1}$, then $m\geq10$ and
\begin{equation*}
\begin{split}
&\sum\nolimits_{i=\min\Gamma_{m,1}}^nd(i)\\
=&\begin{cases}
\sum_{i=\min\Gamma_{m-3,1}}^{n-t_{m-1}}d(i),
&\hat{\theta}_m^1\leq n<\hat{\eta}_m^1;\\
\sum_{i=\min\Gamma_{m-3,1}}^{n-t_{m-1}}d(i)+n-\tfrac{t_{m}+5t_{m-4}}{2}+1,
&\hat{\eta}_m^1\leq n<\hat{\theta}_m^2;\\
\sum_{i=\min\Gamma_{m-2,1}}^{n-t_{m-1}}d(i)
+\sum d(\Gamma_{m-3,1})+\tfrac{t_{m-2}-3t_{m-4}-1}{2},
&\hat{\theta}_m^2\leq n<\hat{\theta}_m^3;\\
\sum_{i=\min\Gamma_{m-1,1}}^{n-t_{m-1}}d(i)+\sum d(\Gamma_{m-3,1})
+\sum d(\Gamma_{m-2,1})+\tfrac{t_{m-2}-3t_{m-4}-1}{2},
&otherwise.
\end{cases}
\end{split}
\end{equation*}
\end{property}

\begin{aD}[]
Step 1. For $n\leq325$, calculate $\sum_{i=1}^n d(i)$ by Property \ref{d}.

Step 2. For $n\geq326$, find the $m$ such that $n\in\Gamma_{m,1}$, then $m\geq10$.
Calculate $\sum_{i=\min{\Gamma_{m,1}}}^n d(i)$ by Property \ref{P6.3}.
And calculate $D(\min{\Gamma_{m,1}}-1)=D(\max{\Gamma_{m-1,1}})$ by Property \ref{P6.2}.

Step 3. $D(n)=D(\max{\Gamma_{m-1,1}})+\sum_{i=\min{\Gamma_{m,1}}}^n d(i)$.
\end{aD}

\noindent\emph{Example.} Consider $n=500\in\Gamma_{10,1}=[326,\cdots,599]$.
Since $\hat{\theta}_{10}^3=\tfrac{3t_{9}+t_{7}-1}{2}=451
\leq n<\hat{\theta}_{11}^1=\tfrac{t_{11}+t_{9}-1}{2}=600$,
by Property 6.3,
$\sum_{i=\min{\Gamma_{10,1}}}^{500} d(i)
=\sum_{i=\min{\Gamma_{9,1}}}^{500-t_9}d(i)
+\sum d(\Gamma_{7,1})
+\sum d(\Gamma_{8,1})+\tfrac{t_{8}-3t_{6}-1}{2}$.
Here $\sum_{i=\min{\Gamma_{9,1}}}^{500-t_9}d(i)=\sum_{i=177}^{226}d(i)
=d([197,198,199,200])+d(207)=5$;
$\sum d(\Gamma_{7,1})=1$,
$\sum d(\Gamma_{8,1})=3$ and $\tfrac{t_{8}-3t_{6}-1}{2}=8$.
Thus $\sum_{i=\min{\Gamma_{10,1}}}^{500} d(i)=17$.
By Property 6.2,
$D(\max{\Gamma_{9,1}})=\frac{9}{44}(9t_{9}-12t_{8}-5t_{7})
+\frac{3}{11}(-2t_{9}+2t_{8}+t_{7})+\frac{9}{4}=12$.
Thus $D(500)=D(\max{\Gamma_{9,1}})+\sum_{i=\min{\Gamma_{10,1}}}^{500} d(i)=12+17=29$.

In fact, the positions of the last letters of the 29 cubes are
$\{58,107,108,139,197,198,199,200,\\207,256,257,288,332,
362,363,364,365,366,367,368,369,381,382,413,471,472,473,474,481\}.$

\vspace{0.2cm}

Now we turn to give the expressions of $D(t_m)$.

For $m\geq10$, $\hat{\theta}_m^3\leq t_{m}<\hat{\theta}_{m+1}^1$
and $\hat{\theta}_{m-1}^2\leq t_{m}-t_{m-1}<\hat{\theta}_{m-1}^3$. By Property \ref{P6.3},
$$\begin{array}{rl} &\sum_{i=\min\Gamma_{m,1}}^{t_m}d(i)-\sum_{i=\min\Gamma_{m-3,1}}^{t_{m-3}}d(i)\\
=&\sum d(\Gamma_{m-3,1})
+\sum d(\Gamma_{m-2,1})+\tfrac{t_{m-2}-3t_{m-4}-1}{2}
+\sum d(\Gamma_{m-4,1})+\tfrac{t_{m-3}-3t_{m-5}-1}{2}\\
=&\sum d(\Gamma_{m-1,1})+\tfrac{t_{m-2}-3t_{m-4}-1}{2}\\
=&\frac{m}{22}(7t_{m-1}-13t_{m-2}+t_{m-3})
+\frac{1}{44}(-121t_{m-1}+188t_{m-2}+57t_{m-3})-\frac{1}{4}\\
=&\frac{m}{22}(t_{m}+6t_{m-1}-14t_{m-2})
+\frac{1}{44}(57t_{m}-178t_{m-1}+131t_{m-2})-\frac{1}{4}.
\end{array}$$
For $m\geq7$, by induction, $\sum_{i=\min\Gamma_{m,1}}^{t_m}d(i)$ is equal to
$$\begin{array}{rl}
\sum_{i=\min\Gamma_{m,1}}^{t_m}d(i)
=&\tfrac{m}{44}(-7t_m+2t_{m-1}+21t_{m-2})
+\tfrac{1}{11}(-10t_m+21t_{m-1}-6t_{m-2})-\frac{m}{12}\\
&~~~~~~+0[m\equiv0(\mathrm{mod}~3)]+\frac{1}{3}[m\equiv1(\mathrm{mod}~3)]
+\frac{2}{3}[m\equiv2(\mathrm{mod}~3)].
\end{array}$$
Here $[P]$ is Iverson notation, and equals 1 if P holds and 0 otherwise.

By Property \ref{P6.2} and $D(t_m)=D(\max\Gamma_{m-1,1})+\sum_{i=\min\Gamma_{m,1}}^{t_m}d(i)$, we prove Theorem 21 in H.Mousavi and J.Shallit\cite{MS2014} in a novel way, i.e. \textbf{Theorem D}.

\end{CJK*}
\end{document}